\documentclass[10pt]{amsart}
\usepackage{amsmath}
\usepackage{amssymb}
\usepackage{indentfirst}
\usepackage[latin1]{inputenc}
\usepackage{geometry}
\usepackage{amsfonts}
\DeclareSymbolFont{largesymbols}{OMX}{yhex}{m}{n}
\DeclareMathAccent{\widehat}{\mathord}{largesymbols}{"62}

\newcommand{\Aut}{\mbox{\rm Aut}}

\newcommand{\Ker}{\mbox{\rm Ker }}

\newcommand{\MM}{{\mathcal{M}}}
\newcommand{\U}{{\mathcal{U}}}

\newcommand{\Z}{{\mathbb Z}}

\newcommand{\Q}{{\mathbb Q}}

\newcommand{\C}{{\mathbb C}}

\newcommand{\HQ}{{\mathbb H}}

\newcommand{\bc}{\begin{center}}
\newcommand{\ec}{\end{center}}

\newcommand{\Or}{{\mathcal O}}
\newcommand{\GL}{{\rm GL}}
\newcommand{\SL}{{\rm SL}}

\newcommand{\Cen}{{\rm Cen}}

\newcommand{\tr}{{\rm tr}}

\newcommand{\Irr}{{\rm Irr}}
\newcommand{\Gal}{\rm Gal}
\newcommand{\GEN}[1]{\left\langle #1 \right\rangle}
\newtheorem{theorem}{Theorem}[section]

\newtheorem{lemma}[theorem]{Lemma}
\newtheorem{corollary}[theorem]{Corollary}
\newtheorem{proposition}[theorem]{Proposition}
\newtheorem{remark}[theorem]{Remark}
\newtheorem{notation}[theorem]{Notation}

\newtheorem{examples}[theorem]{Examples}

\begin{document}

\title[Rational group algebras of finite  groups]{Rational group algebras of finite groups:
from idempotents to units of integral group rings}



\thanks{Research partially supported by Onderzoeksraad of
Vrije Universiteit Brussel, Fonds voor Wetenschappelijk Onderzoek (Flanders),
the grant PN-II-ID-PCE-2007-1 project ID\_532, contract no.~29/28.09.2007,
Ministerio de Ciencia y Tecnolog\'ia of Spain and Fundaci\'{o}n S\'{e}neca of
Murcia.}

\author{Eric Jespers}

\address{Department of Mathematics, Vrije Universiteit Brussel,
Pleinlaan 2, 1050 Brussels, Belgium}

\email{efjesper@vub.ac.be}

\author{Gabriela Olteanu}

\address{Department of Statistics-Forecasts-Mathematics, Babe\c s-Bolyai
University, Str. T. Mihali 58-60, 400591 Cluj-Napoca, Romania}

\email{gabriela.olteanu@econ.ubbcluj.ro}

\author{\'{A}ngel del R\'io}

\address{Departamento de Matem\'{a}ticas, Universidad de Murcia,  30100 Murcia, Spain}

\email{adelrio@um.es}

\date{}

\begin{abstract}
We give an explicit and character-free construction of a complete set of
orthogonal primitive idempotents of a rational group algebra of a finite
nilpotent group and a full description of the Wedderburn decomposition of such
algebras. An immediate consequence is a well-known result of Roquette on the
Schur indices of the simple components of group algebras of finite nilpotent
groups. As an application, we obtain that the unit group of the integral group
ring $\Z G$ of a finite nilpotent group $G$ has a subgroup of finite index that
is generated by three nilpotent groups for which we have an explicit
description of their generators. Another application is a new construction of
free subgroups in the unit group. In all the constructions dealt with, pairs of
subgroups $(H,K)$, called strong Shoda pairs, and explicit constructed central
elements $e(G,H,K)$ play a crucial role. For arbitrary finite groups we prove
that the primitive central idempotents of the rational group algebras are
rational linear combinations of such $e(G,H,K)$, with $(H,K)$ strong Shoda
pairs in subgroups of $G$.
\end{abstract}

\maketitle

\section{Introduction}

The investigation of the unit group $\U(\Z G)$ of the integral group ring $\Z
G$ of a finite group $G$ has a long history and goes back to work of Higman
\cite{H} and Brauer \cite{Bra}. One of the reasons for the importance of the
integral group ring $\Z G$ is that it is  an algebraic tool that links group
and ring theory. It was anticipated for a long time that the  defining group
$G$ would be determined by its integral group ring, i.e. if $\Z G$ is
isomorphic with $\Z H$ for some finite group $H$ then $G\cong H$, the
isomorphism problem. Roggenkamp and Scott \cite{RS} showed that this indeed is
the case if $G$ is a nilpotent group. Weiss proved a more general result
\cite{W}, which also confirmed a Zassenhaus conjecture. It was a  surprise when
Hertweck \cite{He} gave a counter example to the isomorphism problem. In all
these investigations the unit group $\U (\Z G)$ of $\Z G$ is of fundamental
importance. There is a vast literature on the topic. For a survey up to 1994,
the reader is referred to the books of Passman and Sehgal \cite{Pas,Seh1,Seh2}.
Amongst many others, during the past 15 years, the following problems have
received a lot of attention  (we include some guiding references): construction
of generators for a subgroup of finite index in $\U (\Z G)$ \cite{RS1,RS2,JL1},
construction of free subgroups \cite{HP,MS2}, structure theorems for $\U (\Z
G)$ for some classes of groups $G$ \cite{JPRRZ}.

Essential in these investigations is to consider $\Z G$ as a $\Z$-order in the
rational group algebra $\Q G$ and to have a detailed understanding of the
Wedderburn decomposition of $\Q G$. To do so, a first important step is to
calculate the primitive central idempotents of $\Q G$. A classical method for
this  is to apply Galois descent on the primitive central idempotents of the
complex group algebra $\C G$. The latter idempotents are the elements of the
form $e(\chi )=\frac{\chi(1)}{|G|}\sum_{g\in G}\chi (g^{-1}) g$, where $\chi$
runs through the irreducible (complex) characters $\chi$ of $G$. Hence the
primitive central idempotents of $\Q G$ are the elements of the form
$\sum_{\sigma \in \text{Gal}(\Q (\chi )/\Q )} e(\sigma\circ\chi)$ (see for
example \cite{Yam}). Rather recently, Olivieri, del~R\'io and Sim{\'{o}}n
\cite{ORS} obtained a character free method to describe the primitive central
idempotents of $\Q G$ provided $G$ is a monomial group, that is, every
irreducible character of $G$ is induced from a linear character of a subgroup
of $G$. The new method  relies on a theorem of Shoda on pairs of subgroups
$(H,K)$ of $G$ with $K$ normal in $H$, $H/K$ abelian and so that an irreducible
linear character of $H$ with kernel $K$ induces an irreducible character of
$G$. Such pairs are called Shoda pairs of $G$. In Section~2 we recall the
necessary background and explain the description of the primitive central
idempotents of $\Q G$. It turns out that these idempotents can be built using
the central elements $e(G,H,K)$ (see Section~2 for the definition) with $(H,K)$
a Shoda pair of $G$. In case the Shoda pair satisfies some extra conditions
(one calls such a pair a strong Shoda pair) then one also obtains a detailed
description of the Wedderburn component associated with the central idempotent.
This is an important second step towards a description of the simple components
of $\Q G$. This method is applicable to all abelian-by-supersolvable finite
groups, in particular to finite nilpotent groups.

For arbitrary finite groups $G$, it remains an open problem to give a character
free description of the primitive central idempotents of $\Q G$. Only for very
few groups that are not monomial, such a description has been obtained (see for
example \cite{GJ} for alternating groups). In section~3 we show that for
arbitrary finite groups $G$ the elements $e(G,H,K)$ are building blocks for the
construction of the primitive central idempotents $e$ of $\Q G$, i.e. every
such $e$ is a rational linear combination of $e(G,H,K)$, where $(H,K)$ runs
through strong Shoda pairs in subgroups of $G$. The proof makes fundamental use
of Brauer's Theorem on Induced Characters. Presently we are unable to control
the rational coefficients in this linear combination.

In case $G$ is an abelian-by-supersolvable finite group, then, as mentioned
above, the primitive central idempotents of $\Q G$ are of the form $e(G,H,K)$,
with $(H,K)$ a strong Shoda pair of $G$ and the simple component $\Q G
e(G,H,K)$ is described. For nilpotent groups $G$ we will show to have a much
better and detailed control. Indeed, in Section~4 we  describe a complete set
of matrix units (in particular, a complete set of orthogonal primitive
idempotents) of $\Q Ge(G,H,K)$; a third step in the description of $\Q G$. This
allows us to give concrete representations of the projections $ge$ of the group
elements $g\in G$ as matrices over division rings. As a consequence, the
recognition of the $\Z$-order $\Z G$ in the Wedderburn description of $\Q G$ is
reduced to a linear algebra problem over the integers. We include some examples
to show that the method can not be extended to, for example, finite metacyclic
groups. It remains a challenge to construct a complete set of primitive
idempotents for such groups.

In Section~5,  we give several applications to the unit group $\U (\Z G)$ for
$G$ a finite nilpotent group. First we show that if $G$ is a finite nilpotent
group such that $\Q G$ has no exceptional components (see Section~5 for the
definition) then $\U (\Z G)$ has a subgroup of finite index that is generated
by three nilpotent finitely generated groups of which we give explicit
generators. The problem of describing explicitly a finite set of generators for
a subgroup of finite index in $\U (\Z G)$ has been investigated in a long
series of papers. Bass and Milnor did this for abelian groups \cite{Bas}, the
case of nilpotent groups so that their rational group algebra has no
exceptional components was done by Ritter and Sehgal \cite{RS1,RS2}, arbitrary
finite groups  so that their rational group algebra has no exceptional
components were dealt with by Jespers-Leal \cite{JL1}. It was shown that the
Bass cyclic units together with the bicyclic units generate a subgroup of
finite index. Some cases with exceptional components have also been considered,
see for example \cite{GS,J,JL2,Seh3}). In general, very little is known on the
structure of the group generated by the Bass cyclic units and the bicyclic
units, except that ``often'' two of them generate a free group of rank two (see
for example \cite{GP,GR,J,JRR,MS2,JL3}). In this paper, for $G$ a finite
nilpotent group, we not only give new generators for a subgroup of finite
index, but more importantly, the generating set is divided into three subsets,
one of them generating a subgroup of finite index in the central units and each
of the other two generates a nilpotent group. One other advantage of our method
with respect to the proofs and results given in \cite{JL1,RS1,RS2} is that our
proofs are (modulo the central units) more direct and constructive to obtain an
explicit set of generators for a subgroup of finite index in the group of units
of the integral group ring of a finite nilpotent group. Furthermore, we also
give new explicit constructions of free subgroups of rank two.

\section{Preliminaries}
We introduce some useful notation and results, mainly from \cite{JLP} and
\cite{ORS}. Throughout $G$ is a finite group. If $H$ is a subgroup of a group
$G$, then let $\widehat{H}=\frac{1}{|H| }\sum_{h\in H} h\in \Q G$. For $g\in
G$, let $\widehat{g}=\widehat{\langle g \rangle}$ and for non-trivial $G$, let
$\varepsilon(G)=\prod (1-\widehat{M})$, where $M$ runs through the set of all
minimal normal nontrivial subgroups of $G$. Clearly, $\widehat{H}$ is an
idempotent of $\Q G$ which is central if and only if $H$ is normal in $G$. If
$K\lhd H\leq G$ then let \vspace{-1mm}$$\varepsilon(H,K)=\prod_{M/K\in\MM(H/K)}
(\widehat{K}-\widehat{M}),$$ where $\MM(H/K)$ denotes the set of all minimal
normal subgroups of $H/K$. We extend this notation by setting
$\varepsilon(K,K)=\widehat{K}$. Clearly $\varepsilon(H,K)$ is an idempotent of
the group algebra $\Q G$. Let $e(G,H,K)$ be the sum of the distinct
$G$-conjugates of $\varepsilon(H,K)$, that is, if $T$ is a right transversal of
${\rm Cen}_G(\varepsilon(H,K))$ in $G$, then
    $$e(G,H,K)=\displaystyle\sum_{t\in T}\varepsilon(H,K)^t,$$
where $\alpha^{g}=g^{-1}\alpha g$ for $\alpha \in
\C G$ and $g\in G$. Clearly, $e(G,H,K)$ is a
central element of $\mathbb{Q}G$. If the
$G$-conjugates of $\varepsilon(H,K)$ are
orthogonal, then $e(G,H,K)$ is a central
idempotent of $\mathbb{Q}G$.

A Shoda pair of a finite group $G$ is a pair $(H,K)$ of subgroups of $G$ with
the properties that $K\unlhd H$, $H/K$ is cyclic, and if $g\in G$ and
$[H,g]\cap H\subseteq K$ then $g\in H$. A strong Shoda pair of $G$ is a Shoda
pair $(H,K)$ of $G$ such that $H\unlhd N_G(K)$ and the different conjugates of
$\varepsilon(H,K)$ are orthogonal. We also have, in this case, that ${\rm
Cen}_G(\varepsilon(H,K))=N_G(K)$ and $H/K$ is a maximal abelian subgroup of
$N_G(K)/K$ \cite{ORS}.

If $\chi$ is a monomial character of $G$ then $\chi=\psi^G$, the induced
character of a linear character $\psi$ of a subgroup $H$ of $G$. By a Theorem
of Shoda, a monomial character $\chi=\psi^G$ as above is irreducible if and
only if $(H,\Ker \psi)$ is a Shoda pair (see \cite{Shoda} or \cite[Corollary
45.4]{CR}). A finite group $G$ is monomial if every irreducible character of
$G$ is monomial and it is strongly monomial if every irreducible character of
$G$ is strongly monomial. It is well known that every abelian-by-supersolvable
group is monomial (see \cite[Theorem 24.3]{Hup}) and in \cite{ORS} it is proved
that it is even strongly monomial. We will use these results in order to study
the primitive idempotents of group algebras for some abelian-by-supersolvable
groups, including finite nilpotent groups. We will also use the following
description of the simple component associated to a strong Shoda pair.

\begin{theorem}\label{simplessp}\cite[Proposition 3.4]{ORS}
If $(H,K)$ is a strong Shoda pair of G then
$$\Q Ge(G,H,K)\cong M_r(\Q N\varepsilon(H,K))\cong M_r(\Q(\zeta_m)*^{\alpha}_{\tau} N/H),$$
where
   $$m=[H:K], \quad N=N_G(K), \quad r=[G:N]$$
and the action $\alpha$ and twisting $\tau$ are given by
   $$\alpha(nH)(\zeta_m)=\zeta_m^i, \quad \tau(nH,n'H)=\zeta_m^j,$$
if $n^{-1}hnK=h^iK$ and $[n,n']K=h^j K$, for $hK$ a generator of $H/K$,
$n,n'\in N$ and $i,j\in \Z$.
\end{theorem}

In the above theorem, $\zeta_m$ denotes a primitive $m$-th root of unity and we
have used the notation $L*^{\alpha}_{\tau} G$, for $L$ a field and $G$ a group,
to denote a crossed product with action $\alpha:G\rightarrow \Aut(L)$ and
twisting $\tau:G\times G \rightarrow L^*$ \cite{Pas}, i.e. $L*^{\alpha}_{\tau}
G$ is the associative ring $\bigoplus_{g\in G} L u_g$ with multiplication given
by the following rules:
    $$u_g a = \alpha_g(a) u_g, \quad u_g u_h = \tau(g,h) u_{gh}.$$
If the action of $G$ on $L$ is faithful then one may identify $G$ with a group
of automorphisms of $L$ and the center of $L*^{\alpha}_{\tau} G$ is the fixed
subfield $F=L^G$, so that $G=\Gal(L/F)$, and this crossed product is usually
denoted by $(L/F,\tau)$ \cite{Rei}. We refer to these crossed products as
classical crossed products. This is the case for the crossed product $\Q
N\varepsilon(H,K))\cong \Q(\zeta_m)*^{\alpha}_{\tau} N/H$ in
Theorem~\ref{simplessp} which can be described as $(\Q(\zeta_m)/F,\tau)$, where
$F$ is the center of the algebra, which is determined by the Galois action
given in Theorem~\ref{simplessp}.

\section{Primitive central idempotents}

For an irreducible character $\chi$ of $G$ and a field $F$ of characteristic
$0$, $e_F(\chi)$ denotes the only primitive central idempotent of $FG$ such
that $\chi(e)\ne 0$. In this section, using Brauer's Theorem on Induced
Characters, we give a description of every primitive central idempotent
$e_{\Q}(\chi)$ of a rational group algebra $\Q G$ corresponding to an
irreducible character $\chi$ of a finite group $G$ as a rational linear
combination of elements of the form $e(G,H_i,K_i)$, with each  $(H_i,K_i)$ a
strong Shoda pair in a subgroup of $G$,  or equivalently, $K_i$ is a normal
subgroup of $H_i$ with $H_i/K_i$ cyclic.

\begin{theorem}[Brauer]\cite{Bra} Every complex character $\chi$ of a finite group $G$ is a
$\Z$--linear combination $\chi=\sum_i a_i \theta_i^G$ , $a_i\in \Z$, of
characters induced from linear characters $\theta_i$ of elementary subgroups
$M_i$ of $G$, where by an elementary subgroup of $G$ we mean one which is a
direct product of a cyclic group and a $p$-group for some prime $p$.
\end{theorem}

In particular, the $M_i$'s are cyclic--by--$p_i$-groups for some primes $p_i$,
hence by \cite{ORS} each $M_i$ is strongly monomial. As a consequence, every
irreducible character of such a subgroup $M_i$ is an induced character
$\theta_i^{M_i}$ from a linear character $\theta_i$ of a subgroup $H_i$ of
$M_i$. So, $\theta_i^{M_i}$ is irreducible and
$(H_i,\ker(\theta_i))$ is a strong Shoda pair of $M_i$.

\medskip We also will use the result \cite[Theorem~2.1.]{ORS} that describes the
primitive central idempotents $e_{\Q}(\psi^G)$ of a rational group algebra $\Q
G$ associated to a monomial irreducible character $\psi^G$  as
\begin{eqnarray} \label{eQ}
e_{\Q}(\psi^G)=\frac{[\Cen_G(\varepsilon(H,K)):H]}{[\Q(\psi):\Q(\psi^G)]}e(G,H,K)\end{eqnarray}
where $\psi$ is a linear character of the subgroup $H$ of $G$ and $K$ is the
kernel of $\psi$.

\begin{proposition} \label{formulaPCI}
Let $G$ be a finite group of order $n$ and $\chi$ an irreducible character of
$G$. Then the primitive central idempotent $e_{\Q}(\chi)$ of $\Q G$ associated
to $\chi$ is of the form
$$e_{\Q}(\chi)=\frac{1}{[\Q(\zeta_n):\Q(\chi)]}\sum_i
a_i[\Q(\zeta_n):\Q(\psi_i)][C_i:H_i]
e(G,H_i,K_i)$$ where $a_i\in \Z$, $(H_i,K_i)$ are
strong Shoda pairs of subgroups of $G$
(equivalently $H_i/K_i$ is a cyclic section of
$G$), $C_i=\Cen_G(\varepsilon(H_i,K_i))$ and
$\psi_i$ are linear characters of $H_i$ with
kernel $K_i$.
\end{proposition}

\begin{proof}
As it was mentioned in the introduction, for every $\chi\in \Irr(G)$, we have
$$e_{\Q}(\chi)=\displaystyle\sum_{\sigma\in
\text{Gal}(\Q(\chi)/\Q)}e(\chi^{\sigma})=
\displaystyle\sum_{\sigma\in
\text{Gal}(\Q(\chi)/\Q)}\sigma(e(\chi))=\tr_{\Q(\chi)/\Q}(e(\chi)),$$
where $\chi^{\sigma}$ is the character of $G$
given by $\chi^{\sigma}(g)=\sigma(\chi(g))$, for
$g\in G$. The interpretation of $e_{\Q}(\chi)$ as
a trace, suggests the following useful notation
for the next arguments. For any finite Galois
extension $F$ of $\Q$ containing $\Q(\chi)$, let
$$e_{\Q}^{F}=\displaystyle\sum_{\sigma\in
\text{Gal}(F/\Q)}\sigma(e(\chi))=\tr_{F/\Q}(e(\chi)).$$
Hence $e_{\Q}(\chi)=e_{\Q}^{\Q(\chi)}(\chi) =
\frac{1}{[F:\Q(\chi)]}(e_{\Q}^F(\chi))$ for every
finite Galois extension $F$ of $\Q(\chi)$. Using
Brauer's Theorem on Induced Characters, we now
may write $\chi=\sum_i a_i \psi_i^G$, with
$\psi_i$ linear characters of elementary
subgroups $H_i$ with kernel $K_i$ and $a_i\in
\Z$. Then $$e_{\Q}^{\Q(\zeta_n)}(\chi)=\sum_i a_i
e_{\Q}^{\Q(\zeta_n)}(\psi_i^G)$$ and, for every
i, we will compute
$e_{\Q}^{\Q(\zeta_n)}(\psi_i^G)$, as in the proof
of \cite[Proposition 2.1.]{ORS}. (Note that
$\Q(\zeta_n)$ contains $\Q(\chi)$, because it is
an splitting field of $G$.)

Put $e_i=e(\psi_i)$. We know that $\mathcal{A}=\text{Aut}(\C)$ acts on the left
and $G$ acts on the right on $\psi_i$ and on $e_i$ (by composition and by
conjugation respectively) and that their actions are compatible. Hence one may
consider $\mathcal{A}\times G$ acting on the left on the set of irreducible
characters of subgroups of $G$ (and similarly on the $e_i$'s) by
$(\sigma,g)\cdot\psi_i=\sigma\cdot \psi_{i} \cdot g^{-1}$.

Let $\text{Gal}(\Q(\zeta_n)/\Q) = \{\sigma_1,\dots, \sigma_l\}$ and
$T_i=\{g_1,\dots,g_m\}$ a right transversal of $H_i$ in $G$. Denote by
$C_i=\Cen_G(\varepsilon(H_i,K_i))$. We have that $\sum_{k=1}^m e_i\cdot
g_k=e(\psi_i^G)$, hence
\begin{eqnarray*}\label{eq1} e_{\Q}^{\Q(\zeta_n)}(\psi_i^G) &=& \sum_{j=1}^l
\sigma_j e(\psi_i^G) = \sum_{j=1}^l \sum_{k=1}^m \sigma_j e_i\cdot g_k =
\sum_{k=1}^m \tr_{\Q(\zeta_n)/\Q}(e_i)\cdot g_k  \\ &=&\sum_{k=1}^m
[\Q(\zeta_n):\Q(\psi_i)]\tr_{\Q(\psi_i)/\Q}(e_i)\cdot g_k = \sum_{k=1}^m
[\Q(\zeta_n):\Q(\psi_i)]\varepsilon(H_i,K_i)^{g_k}   \\ &=&
[\Q(\zeta_n):\Q(\psi_i)] [C_i:H_i] e(G,H_i,K_i)
\end{eqnarray*}
The above computations now easily yield the desired formula for $e_{\Q}(\chi)$.
\end{proof}

\begin{remark} \rm Notice that the formula from Proposition~\ref{formulaPCI}
for the computation of the primitive central idempotents $e_{\Q}(\chi)$ of $\Q
G$ associated to an irreducible character $\chi$ of $G$ coincides with
formula~(\ref{eQ}) in case $\chi$ is a monomial irreducible character of $G$,
that is $\chi$ is induced to $G$ from only one linear character $\psi$ of a
subgroup $H$, with kernel $K$ such that $(H,K)$ is a Shoda pair of $G$.

In general, as seen in Proposition~\ref{formulaPCI}, one has to consider all
strong Shoda pairs in subgroups of $G$ that contribute to the description of a
primitive central idempotent of $\Q G$. However, one can reduce the search of
the Shoda pairs that determine the primitive central idempotents of $\Q G$ to
representatives given by a relation between such pairs of subgroups. Indeed, in
\cite[Proposition~1.4]{ORS-aut}, it is proved that if $(H_1,K_1)$ and
$(H_2,K_2)$ are two Shoda pairs of a finite group $G$ and $\alpha_1,\alpha_2\in
\Q$ are such that $e_i=\alpha_i e(G,H_i,K_i)$ is a primitive central idempotent
of $\Q G$ for $i=1,2$, then $e_1=e_2$ if and only if there is $g\in G$ such
that $H_1^g\cap K_2=K_1^g\cap H_2$.
\end{remark}

\begin{remark}
\rm We would like to be able to give a bound for the integers $a_i$ used in the previous proposition and one would also
like to give more information on the pairs of groups $(H_i,K_i)$ that one has to consider in the description of
$e_{\Q}(\chi)$.

Notice that for monomial (respectively strongly monomial) groups, all primitive central idempotents are realized as
elements of the form $\alpha e(G,H,K)$, with $\alpha\in \Q$, for some Shoda pair $(H,K)$ (respectively strong Shoda
pair and $\alpha=1$) in $G$. However, for the smallest non-monomial group, which is $\SL(2,3)$, this is not true any
more. Indeed, in \cite[Example 5.7.]{ORS}, the two primitive central idempotents corresponding to the non-monomial
characters are $e_1=\frac{1}{2}e(G,B,A)$, and $e_2=\frac{1}{4}e(G,B,1)-\frac{1}{4}e(G,B,A)$, with $G=\langle
x,y\rangle\rtimes \langle a \rangle$,   a semidirect product of the quaternion group $\langle x,y\rangle$ of order $8$
by the cyclic group $A=\langle a \rangle$ of order $3$, and with  $B=\langle x^2a\rangle$. However, $e_{2}$ can not be
written as a rational linear multiple of some $e(G,H,K)$ with $(H,K)$ a pair of subgroups of $G$ such that $K\unlhd H$.
\end{remark}

\section{Primitive idempotents for finite nilpotent groups}

We start this section by showing a method to produce a complete set of
orthogonal primitive idempotents of a classical crossed product with trivial
twisting $\tau=1$, i.e. $\tau(g,h)=1$, for every $g,h\in G$. Let $L$ be a field
of characteristic zero. Observe that $(L/F,1)\simeq M_n(F)$, with $n=[L:F]$,
therefore a complete set of orthogonal primitive idempotents of $(L/F,1)$
contains $n$ idempotents.

\begin{lemma}\label{TrivialTwist}
Let $A=(L/F,1)$ be a classical crossed product with trivial twisting and let $G=\Gal(L/F)$ with $n=|G| $. Let
$e=\frac{1}{|G| }\sum_{g\in G} u_g$ and let $x_1,\dots,x_n$ be non-zero elements of $L$. Then the conjugates of $e$ by
$x_1,\dots,x_n$ form a complete set of orthogonal primitive idempotents of $A$ if and only if
$\tr_{L/F}(x_ix_j^{-1})=0$ for every $i\ne j$. ($\tr_{L/F}$ is the trace of $L$ over $F$.)
\end{lemma}

\begin{proof}
As the twisting is trivial, $\{u_g:g\in G\}$ is a subgroup of order $|G| $ of the group of units of $A$ and hence $e$
is an idempotent of $A$. Moreover $u_ge=e$ for every $g\in G$. Therefore, if $x \in L$ then $exe= \frac{1}{|G|
}\sum_{g\in G} u_g x e = \frac{1}{|G| }\sum_{g\in G} x^{g^{-1}} u_g e = \frac{1}{|G| }\sum_{g\in G} x^g e =
\frac{1}{|G| }\tr_{L/F}(x) e$. Thus, if $x\in L$ then $e$ and $xex^{-1}$ are orthogonal if and only if $\tr_{L/F}(x)=0$
and the lemma follows.
\end{proof}

\begin{examples} \label{ExTT} {\rm
\begin{enumerate}
\item In the proof of Theorem~\ref{main} we will encounter some examples of
classical crossed products with trivial twisting with a list $x_1,\dots,x_n$
satisfying the conditions of the previous lemma.

\item Another situation where one can find always such elements correspond to
the case when $L/F$ is a cyclic extension of order $n$ and $F$ contains an
$n$-root of unity. Then $L$ is the splitting field over $F$ of an irreducible
polynomial of $F[X]$ of the form $X^n-a$. If $\alpha\in L$ with $\alpha^n=a$
then $x_1=1,x_2=\alpha,\dots,x_n=\alpha^{n-1}$ satisfy the conditions of
Lemma~\ref{TrivialTwist}. Indeed, the minimal polynomial of $\alpha^i$ over $F$
for $1\le i < n$ is of the form $X^{n/d}-a^{i/d}$ for $d=\gcd(n,i)$ and
therefore $\tr_{L/F}(\alpha^i)=[L:F(\alpha^i)]\tr_{F(\alpha^i)/F}(\alpha^i)=0$
and similarly $\tr_{L/F}(\alpha^{-i})=0$.

\item We now construct an example where there are no elements $x_1,\dots,x_n$
satisfying the conditions of Lemma~\ref{TrivialTwist}. Consider the trivial
cyclic algebra $(L=\Q(\zeta_7)/F=\Q(\sqrt{-7}),1)$ of degree 3. If
$x_1,x_2,x_3$ satisfy the conditions of Lemma~\ref{TrivialTwist} then
$\alpha=x_2x_1^{-1}$ and $\alpha^{-1}=x_1x_2^{-1}$ are non-zero elements of $L$
with zero trace over $F$. This implies that the minimal polynomial of $\alpha$
over $F$ is of the form $X^3-a$ for some $a\in F$. But this implies that $F$
contains  a third root of unity, a contradiction.
\end{enumerate}}
\end{examples}

The groups listed in the following lemma will be the building blocks in the
proof of Theorem~\ref{main}. For $n$ and $p$ integers with $p$ prime, we use
$v_p(n)$ to denote the valuation at $p$ of $n$, i.e. $p^n$ is the maximum
$p$-th power dividing $n$.

\begin{lemma} \label{presentationG}
Let $G$ be a finite $p$-group  which has a maximal abelian subgroup which is
cyclic and normal in $G$. Then $G$ is isomorphic to one of the groups given by
the following presentations:
\begin{center}
\begin{tabular}{l}
$P_1=\GEN{a,b\mid a^{p^n}=b^{p^k}=1,\; b^{-1} a b
= a^r}$, with either $v_p(r-1)=n-k$ or $p=2$ and
$r\not\equiv 1 \mod 4$. \\ $P_2=\GEN{a,b,c\mid
a^{2^n}=1,\; b^{2^k}=1,\; c^2 =1,\; bc=cb,\;
b^{-1} a b = a^r,\; c^{-1}ac=a^{-1}}$, with
$r\equiv 1 \mod 4$.
\\ $P_3=\GEN{a,b,c\mid a^{2^n}=1,\; b^{2^k}=1,\; c^2
=a^{2^{n-1}},\; bc=cb,\; b^{-1} a b = a^r,\;
c^{-1}ac=a^{-1}}$, with $r\equiv 1 \mod 4$.
\end{tabular}
\end{center}
\end{lemma}

Note that if $k=0$ (equivalently, if $b=1$) then the first case correspond to
the case when $G$ is abelian (and hence cyclic), the second case coincides with
the first case with $p=2$, $k=1$ and $r=-1$, and the third case is the
quaternion group of order $2^{n+1}$.

\begin{proof}
Let $A$ be a maximal abelian subgroup of $G$ and assume that $A$ is cyclic
(generated by $a$) and normal in $G$. Put  $|A| =p^n$. Consider the action of
$G$ on $A$ by inner automorphisms. Since $A$ is maximal abelian in $G$, the
kernel of this action is $A$ and therefore $G/A$ is isomorphic to a subgroup of
the group of automorphisms of $A$.

If either $p$ is odd or $p=2$ and $n\le 2$ then $\Aut(A)$ is cyclic and
otherwise $\Aut(A)=\GEN{\phi_5}\times \GEN{\phi_{-1}}$, where $\phi_r$ is the
automorphism of $A$ given by $\phi_r(x)=x^r$.

Assume that $G/A$ is cyclic, so that $G$ has a presentation of the form
    \begin{equation}\label{Meta}
    G=\GEN{a,b\mid a^{p^n}=1, b^{p^k}=a^s,b^{-1} a b=a^r},
    \end{equation}
with $p^n\mid r^{p^k}-1$ and $p^n\mid s(r-1)$. If $p^i\ge 3$ then $(1+xp^i)^p
\equiv 1+xp^{i+1} \mod p^{i+2}$ for every $i\ge 1$ and $x\in \Z$. Using this,
one deduces that if either $p$ is odd or $p=2$ and $r\equiv 1 \mod 4$ then
$v_p(r^{p^i}-1)=v_p(r^{p^{i-1}}-1)+1$, for every $i\ge 1$. Furthermore, from
the assumption that $A$ is maximal abelian in $G$, one deduces that $n\le
v_p(r^{p^k}-1) = v_p(r^{p^{k-1}}-1)+1\le n$ and hence $v_p(r^{p^k}-1)=n$ and
$v_p(r-1)=n-k$. Therefore, $v_p(s)\ge
n-v_p(r-1)=k=v_p\left(\frac{r^{p^k}-1}{r-1}\right) =
v_p\left(1+r+r^2+\dots+r^{p^k-1}\right)$ and hence there is an integer $x$ such
that $x(1+r+r^2+\dots+r^{p^k-1})+s\equiv 1 \mod p^n$. Then $(a^xb)^{p^k}=1$
and, replacing $b$ by $a^xb$ in (\ref{Meta}), we obtain the presentation of
$P_1$. We have also proved that $v_p(r-1)=n-k$ unless $p=2$ and $r\not\equiv 1
\mod 4$. Assume $p=2$ and $r\not\equiv 1 \mod 4$, $v_2(s)\ge n-v_2(r-1)=n-1$
and $v_2(1+r+r^2+\dots r^{2^k-1}) = v_2\left(\frac{r^{2^k}-1}{r-1}\right) \ge
n-1$. If $v_2(s)\ge n$ then $G\cong P_1$. If $v_2(s)=v_2(1+r+r^2+\dots
r^{2^k-1})=n-1$ then $(ab)^{2^{k}}=1$. Replacing $b$ by $ab$ we obtain again
that $G\cong P_1$. Otherwise, $v_2(s)=n-1$ and $v_2(1+r+r^2+\dots
r^{2^k-1})=n$. Therefore $v_2(r^{2^k-1})>n>v_2(r^{2^{k-1}}-1)$ and by the first
part of the proof (applied to $\GEN{a,b^2}$) this implies that $k=2$. Then $G$
is the quaternion group of order $2^{n+1}$ which is $P_3$ for $k=0$.

Assume now that $G/A$ is not cyclic, so $p=2$ and $G/\langle a\rangle=\langle
\overline{b} \rangle \times \langle \overline{c} \rangle$ with $c$ acting by
inversion on $\langle a \rangle$. This provides a presentation of $G$ of the
form
\begin{equation}G=\langle a,b,c \mid  a^{2^n}=1, b^{2^k}=a^s, ca=a^{-1}c, cb=a^ibc,
b^{-1}ab=a^r, c^2=1 \text{ or } c^2=a^{2^{n-1}}\rangle.\end{equation}

Replacing $b$ by $bc$ if needed, one may assume that $v_2(r-1)=n-k\ge 2$ and
$v_2(r^{2^k}-1)=n$. So, applying the first part of the proof to $\GEN{a,b}$, we
may assume that $b^{2^k}=1$. Then
$c=b^{2^k}cb^{-2^k}=a^{-i(1+r+\ldots+r^{2^k-1})}c$ and so $2^n\mid
-i(1+r+\ldots+r^{2^k-1})=-i\frac{r^{2^k}-1}{r-1}$. As
$v_2\left(\frac{r^{2^k}-1}{r-1}\right)=k$, we have $v_2(i)\ge n-k=v_2(r-1)$.
Hence, there exists an integer $j$ so that $j(r-1) -i \equiv 0 \mod 2^n$. It is
easy to verify that the commutator of $b$ and $a^jc$ is 1. So, replacing $c$ by
$a^jc$ if needed, we may assume that $b$ and $c$ commute and we obtain the
presentation of $P_2$, if $c^2=1$, and the presentation of $P_3$, if
$c^2=a^{2^{n-1}}$.
\end{proof}

We also need the following result on splitting of a Hamiltonian quaternion
algebra $\HQ(F)=F[i,j\mid i^2=j^2=-1,ji+ji=0]$.

\begin{lemma}\label{split quaternion}
Let $F$ be a field of characteristic zero. Then the quaternion algebra $\HQ(F)$ splits
if and only if $x^2+y^2=-1$ for some $x,y\in F$. In that case $\frac{1}{2}(1+xi+yj)$ and
$\frac{1}{2}(1-xi-yj)$ form a complete set of primitive idempotents of $\HQ(F)$.

Furthermore, if $F=\Q(\zeta_m,\zeta_{2^n}+\zeta_{2^n}^{-1})$ with $m$ odd then
$-1$ is the sum of two squares of $F$ if and only if $m\ne 1$ and either $n\ge
3$ or the multiplicative order of $2$ modulo $m$ is even.
\end{lemma}

\begin{proof}
The first part can be found in \cite[Proposition 1.13]{Seh1}. Now assume that
$F=\Q(\zeta_m,\zeta_{2^n}+\zeta_{2^n}^{-1})$ with $m$ odd. If $m=1$ then $F$ is
totally real and therefore $-1$ is not the sum of two squares of $F$. So assume
that $m\ne 1$. If $n\le 2$, then $F=\Q(\zeta_m)$ and the result is well known
(see for example \cite{Mos,FGS} or \cite[pages 307--308]{Lam}). Finally assume
that $m\ne 1$ and $n\ge 3$. Then $F$ contains $\sqrt{2}$ and, as $2$ is not a
square in $\Q_2$, the duadic completion of $\Q$ \cite[Corollary 2.24]{Lam}, we
deduce that $[F_2:\Q_2]$ is even. Then $-1$ is a sum of squares in $F_p$ for
every place $p$ of $F$ and hence $-1$ is a sum of squares in $F$ (see
\cite[page 304]{Lam}).
\end{proof}

Now we are ready to show an effective method to calculate a complete set of
orthogonal primitive idempotents of $\Q G$ for $G$ a finite nilpotent group.
Since $G$ is abelian-by-supersolvable and hence strongly monomial, it follows
from \cite[Theorem 4.4]{ORS} that every primitive central idempotent of $\Q G$
is of the form $e(G,H,K)$ with $(H,K)$ a strong Shoda pair of $G$ and therefore
it is enough to obtain a complete set of orthogonal primitive idempotents of
$\Q Ge(G,H,K)$ for every strong Shoda pair $(H,K)$ of $G$. This is described in
our main result that we state now.

\begin{theorem} \label{main} Let $G$ be a finite nilpotent group and $(H,K)$ a
strong Shoda pair of $G$. Set $e=e(G,H,K)$, $\varepsilon=\varepsilon(H,K)$,
$H/K=\langle \overline{a}\rangle$, $N=N_G(K)$ and let $N_2/K$ and
$H_2/K=\langle \overline{a_2}\rangle$ (respectively $N_{2'}/K$ and
$H_{2'}/K=\langle \overline{a_{2'}}\rangle$) denote the $2$-parts
(respectively, $2'$-parts) of $N/K$ and $H/K$ respectively. Then $\langle
\overline{a_{2'}}\rangle$ has a cyclic complement $\langle
\overline{b_{2'}}\rangle$ in $N_{2'}/K$.

A complete set of orthogonal primitive idempotents of $\Q Ge$ consists of the conjugates of
$\widehat{b_{2'}}\beta_2\varepsilon$ by the elements of $T_{2'}T_2T_{G/N}$, where
$T_{2'}=\{1,a_{2'},a_{2'}^2,\dots,a_{2'}^{[N_{2'}:H_{2'}]-1}\}$, $T_{G/N}$ denotes a left transversal of $N$ in $G$ and
$\beta_2$ and $T_2$ are given according to the cases below.

\begin{enumerate} \item  If $H_2/K$ has a complement $M_2/K$ in $N_2/K$
then $\beta_2=\widehat{M_2}$. Moreover, if $M_2/K$ is cyclic then
 there exists $b_2\in N_2$ such that $N_2/K$ is given
by the following presentation
    $$\left\langle \overline{a_2}, \overline{b_2} \mid
    \overline{a_2}^{2^n}=\overline{b_2}^{2^k}=1,\;
    \overline{a_2}^{\overline{b_2}} = \overline{a_2}^r \right\rangle,$$
and if $M_2/K$ is not cyclic, there exist $b_2,c_2\in N_2$ such that $N_2/K$ is given by the following presentation
    $$\left\langle \overline{a_2}, \overline{b_2}, \overline{c_2}  \mid
    \overline{a_2}^{2^n}=\overline{b_2}^{2^k}=1,\;
    \overline{c_2}^2=1,\;
    \overline{a_2}^{\overline{b_2}}=\overline{a_2}^r,\;
    \overline{a_2}^{\overline{c_2}}=\overline{a_2}^{-1},\;
    [\overline{b_2},\overline{c_2}]=1\right\rangle,$$
with $r\equiv 1 \mod 4$ (or equivalently, $\overline{a_2}^{2^{n-2}}$ is central in $N_2/K$). Then
\begin{itemize}
\item[(i)] $T_2=\{1, a_2,a_2^2,\dots,a_2^{2^k-1}\}$, if $\overline{a_2}^{2^{n-2}}$ is central
in $N_2/K$ 
and $M_2/K$ is cyclic; and

\item[(ii)]
$T_2=\{1,a_2,a_2^2,\dots,a_2^{2^{k-1}-1},a_2^{2^{n-2}},a_2^{2^{n-2}+1},\dots,a_2^{2^{n-2}+2^{k-1}-1}\}$,
otherwise.
\end{itemize}

\item if $H_2/K$ has no complement in $N_2/K$ then
there exist $b_2, c_2\in N_2$ such that $N_2/K$ is given by the following presentation
    $$\left\langle \overline{a_2}, \overline{b_2}, \overline{c_2}  \mid
    \overline{a_2}^{2^n}=\overline{b_2}^{2^k}=1,\;
    \overline{c_2}^2=\overline{a_2}^{2^{n-1}},\;
    \overline{a_2}^{\overline{b_2}}=\overline{a_2}^r,\;
    \overline{a_2}^{\overline{c_2}}=\overline{a_2}^{-1},\; [\overline{b_2},\overline{c_2}]=1\right\rangle,$$
with $r\equiv 1 \mod 4$ and we set $m=[H_{2'}:K]/[N_{2'}:H_{2'}]$. Then
\begin{itemize}
\item[(i)] $\beta_2=\widehat{b_2}$ and $T_2=\{1,a_2,a_2^2,\dots,a_2^{2^k-1}\}$,
if either $H_{2'}=K$ or the order of $2$ modulo $m$ is odd and $n-k\le 2$ and

\item[(ii)] $\beta_2=\widehat{b_2}\frac{1+xa_2^{2^{n-2}}+ya_2^{2^{n-2}}c_2}{2}$
and
$T_2=\{1,a_2,a_2^2,\dots,a_2^{2^k-1},c_2,a_2c_2,a_2^2c_2,\dots,a_2^{2^k-1}c_2\}$
with
    $$x,y\in \Q\left[a_{2'}^{\;[N_{2'}:H_{2'}]},a_2^{\;2^k}+a_2^{\;-2^k}\right],$$
satisfying $(1+x^2+y^2)\varepsilon=0$, if $H_{2'}\ne K$ and either the order of $2$ modulo $m$ is
even or $n-k>2$.
\end{itemize}
\end{enumerate}
\end{theorem}
\begin{proof}
We start the proof by making some useful reductions. Taking $T=T_{G/N}$ a left
transversal of $N$ in $G$, the conjugates of $\varepsilon$ by elements of $T$
are the ``diagonal'' elements in the matrix algebra $\Q G e=M_{G/N}(\Q
N\varepsilon)$. Hence, following the proof of \cite[Proposition~3.4]{ORS}, one
can see that it is sufficient to compute a complete set of orthogonal primitive
idempotents for $\Q N\varepsilon=\Q H\,\varepsilon * N/H$ and then add their
$T$-conjugates in order to obtain the primitive idempotents of $\Q G e$. So one
may assume that $N=G$, i.e. $K$ is normal in $G$ and hence $e=\varepsilon$ and
$T=\{1\}$. Then the natural isomorphism $\Q G\widehat{K}\simeq \Q(G/K)$ maps
$\varepsilon$ to $\varepsilon(H/K)$. So, from now on we assume that $K=1$ and
hence $H=\GEN{a}$ is a cyclic maximal abelian subgroup of $G$, which is normal
in $G$ and $e=\varepsilon=\varepsilon(H)$. If $G=H$ then $\Q Ge$ is a field,
$T_2=T_{2'}=\{1\}$ and $b_{2'}=\beta_2=1$; hence the result follows. So, in the
remainder of the proof we assume that $G\ne H$.

The map $a\varepsilon\mapsto \zeta$ induces an isomorphism $f:\Q
H\varepsilon\rightarrow \Q(\zeta)$, where $\zeta$ is a primitive $|H| $-root of
unity. Using the description of $\Q Ge$ given in Theorem~\ref{simplessp}, one
obtains a description of $\Q Ge$ as a classical crossed product
$(\Q(\zeta)/F,\tau)$, where $F$ is the image under $f$ of the center of $\Q
Ge$.

We first consider the case when $G$ is a $p$-group. Then $G$ and $H=\GEN{a}$
satisfy the conditions of Lemma~\ref{presentationG} and therefore $G$ is
isomorphic to one of the three groups of this lemma. Moreover, $H$ has a
complement in $G$ if and only if $G\cong P_1$ or $G\cong P_2$ and, in these
cases, $\tau$ is trivial. We claim that in these cases it is possible to give a
list of elements $x_1,\dots,x_{p^{k}}$ of $\Q(\zeta)$ ($p^{k}=[G:H]$)
satisfying the conditions of Lemma~\ref{TrivialTwist} and the elements
$f^{-1}(x_1),\dots,f^{-1}(x_{p^{k}})$ correspond to the conjugating elements in
$G$ given in the statement of the theorem in the different cases. To prove this
we will use the following fact: if $K$ is a subfield of $\Q(\zeta)$ such that
$\zeta_p\in K$,  $\zeta^i\not\in K$ (with $i=1,\dots,p^k-1$) and, moreover,
$\zeta_4\in K$ if $p=2$   then $\tr_{\Q(\zeta)/K}(\zeta^i)=0$. To see this
notice that if $d$ is the minimum integer such that $\zeta^{ip^d}\in K$ then
$\Q(\zeta^{i})/K$ is cyclic of degree $p^d$ and $\zeta^i$ is a root of
$X^{p^d}-\zeta^{ip^d}\in K[X]$. Then $X^{p^d}-\zeta^{ip^d}$ is the minimal
polynomial of $\zeta^i$ over $K$. Hence $\tr_{\Q(\zeta^{i})/K}(\zeta^i)=0$ and
thus $\tr_{\Q (\zeta)/K} (\zeta^{i})=0$.

 Assume first that $G=P_1$ and $v_p(r-1)=n-k$
(equivalently $a_p^{p^{n-k}}\in Z(G)$), that is, either $p$ is odd or $p=2$ and
$r\equiv 1\mod 4$). Then $F$ is the unique subfield of index $[G:H]=p^k$ in
$\Q(\zeta)$ and such that if $p=2$ then $\zeta_4\in F$. Namely
$F=\Q(\zeta_{p^{n-k}})=\Q(\zeta^{p^k})$. If we set $x_i=\zeta^i$, for
$i=0,1,\dots,p^k-1$, then $x_ix_j^{-1} = \zeta^{i-j}$. If $i\ne j$ then
$\zeta^{i-j}\not\in F$ and hence
$\tr_{\Q(\zeta)/F}(x_ix_j^{-1})=\tr_{\Q(\zeta)/F}(\zeta^{i-j})=0$. Thus, by
Lemma~\ref{TrivialTwist}, the conjugates of $\widehat{b}$ by
$1,\zeta,\zeta^2,\dots,\zeta^{p^k-1}$ form a complete set of orthogonal
primitive idempotents of $(\Q(\zeta)/F,1)$. Then the elements $f^{-1}(x_i)$
form the elements of $T_{2'}$ if $p$ is odd or the elements of $T_2$, in case
(1.i).

Assume now that $G$ is still $P_1$, but with $p=2$ and $r\not\equiv 1 \mod 4$
(equivalently, $a_2^{2^{n-2}}$ is not central). In this case $\zeta_4\not\in F$
and $F(\zeta_4)$ is the unique subextension of $\Q(\zeta)/\Q(\zeta_4)$ of index
$[G:H]/2=2^{k-1}$. That is
$F(\zeta_4)=\Q(\zeta_{2^{n-k+1}})=\Q(\zeta^{2^{k-1}})$. We take $x_i=\zeta^i$
and $x_{2^{k-1}+i}= \zeta^{2^{n-2}+i}=\zeta_4\zeta^i$, for $0\le i < 2^{k-1}$.
Hence, if $i\ne j$ then $x_ix_j^{-1}$ is either $\zeta_4^{\pm 1}$ or
$\zeta^{\pm i}$ or $\zeta_4^{\pm 1} \zeta^{\pm i}$, with
$i=1,2,\dots,2^{k-1}-1$. As $\zeta^i\not\in F(\zeta_4)$, we have
$\tr_{\Q(\zeta)/F(\zeta_4)}(\zeta^i)=0$. Since
$$\tr_{\Q(\zeta)/F}(\zeta_4) =
\tr_{F(\zeta_4)/F}\tr_{\Q(\zeta)/F(\zeta_4)}(\zeta_{4})=[\Q(\zeta):\Q(\zeta_4)]\tr_{F(\zeta_4)/F}(\zeta_4)=0,$$
$$\tr_{\Q(\zeta)/F}(\zeta^i) =
\tr_{F(\zeta_4)/F}\tr_{\Q(\zeta)/F(\zeta_4)}(\zeta^i)=0$$
and $$\tr_{\Q(\zeta)/F}(\zeta_4\zeta^i) =
\tr_{F(\zeta_4)/F}\tr_{\Q(\zeta)/F(\zeta_4)}(\zeta_4
\zeta^i)=\tr_{F(\zeta_4)/F}(\zeta_4\tr_{\Q(\zeta)/F(\zeta_4)}(\zeta^i))
= 0,$$ we deduce that $\tr(x_ix_j^{-1})=0$ for
every $i\ne j$. Then $f^{-1}$ maps these elements
to the elements of $T_2$ for case (1.ii).

Now assume that $G=P_2$. Then $F=\Q(\zeta_{2^{n-k}}+\zeta_{2^{n-k}}^{-1})$.
Since $r\equiv 1 \mod 4$, $n-k\ge 2$. Then the same argument as in the previous
case shows that the $2^{k+1}$ elements of the form $x_i=\zeta^i$ and
$x_{2^k+i}= \zeta^{2^{n/2}+i}=\zeta_4\zeta^i$, for $0\le i < 2^k$ satisfy the
conditions of Lemma~\ref{TrivialTwist}. The elements $f^{-1}(x_i)$ form now the
set $T_2$ of case (1.ii).

Now we consider the non-splitting case, i.e. $G=P_3$. Then the center of $\Q Ge$ is isomorphic to
$F=\Q(\zeta)^{\GEN{b,c}}=\Q(\zeta_{2^{n-k}}+\zeta_{2^{n-k}}^{-1})$ and $\widehat b \,\Q G \,\varepsilon \widehat b =
\widehat b \,\Q \langle a,c \rangle \,\varepsilon \widehat b F+F a^{2^{n-2}}+F c + F(a^{2^{n-2}} c) \cong \HQ(F)$,
which is a division algebra, as $F$ is a real field. Then $\widehat{b}\varepsilon$ is a primitive idempotent of $\Q
Ge$. Hence $\Q Ge \simeq M_{2^k}(\HQ(F))$ and from the first case one can provide the $2^k$ orthogonal primitive
idempotents needed in this case by taking the conjugates of $\widehat{b}$ by $1,a,a^2,\dots,a^{2^k-1}$, and this agrees
with case (2.i). This finishes the $p$-group case.

Let us now consider  the general case, where $G$ is not necessarily a
$p$-group. Then $G=G_2\times G_{p_1}\times \dots \times G_{p_r}=G_2\times
G_{2'}$, with $p_i$ an odd prime for every $i=1,\dots,r$. Then $(H,1)$ is a
strong Shoda pair of $G$ if and only if $(H_{p_i},1)$ is a strong Shoda pair of
$G_{p_i}$, for every $i=0,1,\dots,m$, (with $p_0=2$) and
$\varepsilon(H)=\prod_i\,\varepsilon(H_{p_i})$. Using this and a dimension
argument it easily follows that the simple algebra $\Q G \varepsilon(H)$ is the
tensor product over $\Q$ of the simple algebras $\Q G_{p_i}
\varepsilon(H_{p_i})$. On the other hand, we have seen that for $i\ge 1$, $\Q
G_{p_i} \varepsilon(H_i)\simeq M_{p_i^{k_i}}(\Q(\zeta_{p_i^{n_i-k_i}}))$, for
$p_i^{n_i}=\mid H_i\mid $ and $p_i^{k_i}=[N_i:K_i]$. Then $\Q G_{2'}
\varepsilon(H_{2'}) \simeq M_{[G_{2'}:H_{2'}]}(\Q(\zeta_m))$, with $m=\mid
H_{2'}\mid /[G_{2'}:H_{2'}]$ ($=[H_{2'}:K]/[G_{2'}:H_{2'}]$) and then a
complete set of orthogonal primitive idempotents of $\Q G_{2'}
\varepsilon(H_{2'})$ can be obtained by multiplying the different sets of
idempotents obtained for each tensor factor. Observe that each $G_{p_i}$, with
$i\ge 1$, takes the form $\GEN{a_i}\rtimes \GEN{b_i}$ and so
$G_{2'}=\GEN{a}\rtimes \GEN{b}$, with $a=a_1\dots a_r$ and $b=b_1\dots b_r$.
Having in mind that $a^{p_i^{k_i}}$ is central one can easily deduce, with the
help of the Chinese Remainder Theorem, that the product of the different
primitive idempotents of the factors from the odd part (i.e. the conjugates of
$\widehat{b_i}$ by $1,a_i,a_i^2,\dots,a_i^{p_i^{k_i}-1}$ are the conjugates of
$\widehat{b}\varepsilon$ by $1,a,a^2,\dots,a^{[G_{2'}:H_{2'}]-1}$. In the
notation of the statement of the theorem, $a=a_{2'}$ and
$T_{2'}=\{1,a,a^2,....,a^{[G_{2'}:H_{2'}]-1}$ as wanted.

If $|G| $ is odd then the proof is finished. Otherwise we should combine the
odd and even parts of $G$. If $H_2$ has a complement in $G_2$ then $\Q G_2
\varepsilon(H_2)$ is split over its center and hence we can take $T_2$ as in
the $2$-group case. However, if $H_2$ does not have a complement in $G_2$ then
$\Q G_2
\varepsilon(H_2)=M_{[G_2:H_2]/2}(\HQ(\Q(\zeta_{2^{n-k}}+\zeta_{2^{n-k}}^{-1})))$
and hence $\Q G \varepsilon = M_{[G:H]/2}(\HQ(F))$, with
$F=\Q(\zeta_m,\zeta_{2^{n-k}}+\zeta_{2^{n-k}}^{-1}))$. If $\HQ(F)$ is not split
(equivalently the conditions of (2.i) hold) then we can also take $T_2$ as in
the $2$-group case. However, if $\HQ(F)$ is split then one should duplicate the
number of idempotents, or equivalently duplicate the size of $T_2$. In this
case $-1$ is a sum of squares in $F$. Observing that
$f(a_{2'}^{[N_{2'}:H_{2'}]})$ is a primitive $m$-th root of unity and
$f(a_2^{2^k})$ is a primitive $2^{n-k}$ root of unity, we deduce that there are
$x,y\in \Q(a_{2'}^{[N_{2'}:H_{2'}]},a_2^{2^k}+a_2^{-2^k})$ such that
$(1+x^2+y^2)\varepsilon=0$. Then we can duplicate the number of idempotents by
multiplying the above idempotents by
$f=\frac{1+xa_2^{2^{n-2}}+ya_2^{2^{n-2}}c_2}{2}$ and
$1-f=\frac{1-xa_2^{2^{n-2}}-ya_2^{2^{n-2}}c_2}{2}$ (see Lemma~\ref{split
quaternion}). Observing that $1-f=f^{c_2}$, we obtain that these idempotents
are the conjugates of $\widehat{b_2}f\varepsilon$ by
$1,a_2,\dots,a_2^{2^k-1},c_2,a_2c_2,\dots,a_2^{2^k-1}c_2$, as desired.
\end{proof}

\begin{remark} \label{simple} \rm A description of the simple algebras $\Q Ge$ using
Theorem~\ref{simplessp} can be given according to the cases listed above. Thus,
$\Q Ge=M_{\mid G/H\mid }(\Q(\zeta_{[H:K]})^{N/K})$, that is a matrix algebra
over the fixed field of the natural action of $N/K$ on the cyclotomic field
$\Q(\zeta_{[H:K]})=\Q H\varepsilon(H,K)$, in cases (1) and (2.ii) of
Theorem~\ref{main} and $\Q Ge=M_{\frac{1}{2}\mid G/H\mid }(\HQ
\,(\Q(\zeta_{[H:K]})^{N/K}))$ in case (2.i). In particular, if $\Q Ge$ is a
non-commutative division algebra then $[G:H]=2$, $N=G$ and $\Q Ge\cong
\HQ\left(\Q\left(\zeta_{[H:K]}+\zeta_{[H:K]}^{-1}\right)\right)$, a totally
definite quaternion algebra.
\end{remark}

As a consequence of Theorem~\ref{main}, we get the following result on the
Schur indices of the simple components of group algebras for finite nilpotent
groups over fields of characteristic zero.

\begin{theorem}[Roquette] Let $G$ be a finite nilpotent
group and $F$ a field of characteristic zero. Then $FG\cong \bigoplus_i
M_{n_i}(D_i)$, where $D_i$ are quaternion division algebras if not commutative,
that is the Schur index of the simple components of $FG$ is at most $2$. If the
Schur index of a simple component of $FG$ is $2$ then the Sylow $2$-subgroup of
$G$ has a quaternion section.
\end{theorem}

\begin{remark}
\rm Notice that the use of Lemma~\ref{TrivialTwist} has been essential in all
the cases of the proof of Theorem~\ref{main}. We would like to be able to give
a similar description to the one from Theorem~\ref{main} for a complete set of
orthogonal primitive idempotents for rational group algebras of arbitrary
finite metacyclic groups. Unfortunately, the approach of Theorem~\ref{main}
does not apply here. For example, if $G=C_7\rtimes C_3=\GEN{a}\rtimes \GEN{b}$,
with $b^{-1}ab=a^2$ and $\varepsilon=\varepsilon(\GEN{a})$ then there is not a
complete set of orthogonal primitive idempotents of $\Q G\varepsilon$ formed by
$\Q(a)$-conjugates of $\widehat{b}\varepsilon$. This is a consequence of
Example~\ref{ExTT} (3).
\end{remark}

\begin{notation}\label{Tx}
\rm As an application of Theorem~\ref{main} we next will describe a complete
set of matrix units in a simple component $\Q Ge$, where $e=e(G,H,K)$ with
$(H,K)$ a strong Shoda pair of a finite nilpotent group $G$. A complete set of
primitive idempotents of $\Q Ge(G,H,K)$ is given according to the cases of
Theorem~\ref{main}. Using the notation in these different cases of Theorem 4.5,
let $T_e=T_{2'}T_2T_{G/N}$ and $\beta_e=\widehat{b_{2'}}\beta_2\varepsilon$,
where $\varepsilon=\varepsilon(H,K)$, $T_{G/N}$ denotes a left transversal of
$N=N_G(K)$ in $G$;
$T_{2'}=\left\{1,a_{2'},\dots,a_{2'}^{[N_{2'}:H_{2'}]-1}\right\}$;
    \begin{eqnarray*}\label{T}
    T_2&=&\left\{
    \begin{array}{ll}
        \left\{1, a_2,\dots,a_2^{2^k-1}\right\}, & \hbox{in cases (1.i) and (2.i);} \\
        \left\{1,a_2,\dots,a_2^{2^{k-1}-1},a_2^{2^{n-2}},a_2^{2^{n-2}+1},\dots,a_2^{2^{n-2}+2^{k-1}-1}\right\}, &
        \hbox{in case (1.ii);} \\
        \left\{1,a_2,\dots,a_2^{2^k-1},c_2,a_2c_2,\dots,a_2^{2^k-1}c_2\right\}, & \hbox{in case (2.ii);} \\
    \end{array}
    \right.
     \end{eqnarray*}
and \begin{eqnarray*}\label{x} \beta_2=\left\{%
\begin{array}{ll}
    \widehat{M_2}, & \hbox{in case (1);} \\
    \widehat{b_2}, & \hbox{in case (2.i);} \\
    \widehat{b_2}\frac{1+xa_2^{2^{n-2}}+ya_2^{2^{n-2}}c_2}{2}, &
    \hbox{in case (2.ii).} \\
\end{array}%
\right.\end{eqnarray*}
\end{notation}

\begin{corollary} \label{matrixunits} Let $G$ be a finite nilpotent group.
For every primitive central idempotent
$e=e(G,H,K)$, with $(H,K)$ a strong Shoda pair of
$G$, let $T_e$ and $\beta_e$ be as in
Notation~\ref{Tx}. For every $t,t'\in T_e$ let
$$E_{tt'}= t^{-1}\beta_e \, t'.$$ Then $\{
E_{tt'} \mid  t,t'\in T_e\}$ gives a complete set
of matrix units in $\Q G e$, i.e. $e=\sum_{t\in
T_e} E_{tt}$ and
$E_{t_1t_2}E_{t_3t_4}=\delta_{t_2t_3}E_{t_1t_4}$,
for every $t_1,t_2,t_3,t_4\in T_e$.

Moreover, $E_{tt}\Q G E_{tt}\cong F$, in cases (1) and (2.ii) of
Theorem~\ref{main}, and $E_{tt}\Q G E_{tt}=\HQ(F)$, in case (2.i) of
Theorem~\ref{main}, where $F$ is the fixed subfield of $\Q (a)\varepsilon$
under the natural action of $N/H$.
\end{corollary}

\begin{proof} We know from Theorem~\ref{main} that the set $\{ E_{tt}  \mid   t\in T_e\}$
is a complete set of primitive idempotents of $\Q G e$. From the definition of
the $E_{tt'}$ it easily follows that
$E_{t_1t_2}E_{t_3t_4}=\delta_{t_2t_3}E_{t_1t_4}$, for $t_i\in T_e$, $i=1,\dots,
4$. The second statement is already mentioned in Remark~\ref{simple}.
\end{proof}

\section{Generators of a subgroup of finite index in $\U(\mathbb{Z}G)$}

As an application of the description of the primitive central idempotents in
Corollary~\ref{matrixunits}, we now can easily explicitly construct two
nilpotent subgroups of $\U(\Z G)$ (which correspond with upper respectively
lower triangular matrices in the simple components). Together with the central
units they generate a subgroup of finite index in the unit group. We begin by
recalling a result of Jespers, Parmenter and Sehgal that gives explicit
generators for a subgroup of finite index in the center.

First we recall the construction of units known as Bass cyclic units in the
integral group ring $\Z G$ of a finite group $G$. Let $g\in G$ and suppose $g$
has order $n$. Let $k$ be an integer so that $1 < k < n$ and $(k,n)=1$. Then,
$k^{\varphi(n)}\equiv 1\mod n$, where $\varphi$ is the Euler
$\varphi$-function, and
$$b(g,k)=\left(\sum_{j=0}^{k-1}
g^j\right)^{\varphi(n)}+(1-k^{\varphi(n)})\widehat{g}$$
is a unit in $\Z G$. (Note that our notation
slightly differs from  the one used in
\cite{Seh1}; this because of our definition of
$\widehat{g}=\frac{1}{n}\sum_{i=1}^ng^i$.) The
group generated by all Bass cyclic units of $\Z
G$ we denote by $B(G)$.

We now introduce the following units defined in \cite{JPS}. Let $Z_{v}$ denote
the $v$-th center of $G,$ and suppose that $G$ is nilpotent of class $n$. For
any $x\in G$ and $b\in \Z \langle x\rangle $, put $b_{(1)}=b$ and, for $2\leq
v\leq n$, put
    $$b_{(v)}=\prod_{g\in Z_{v}} b_{(v-1)}^{g}.$$
Note that by induction $b_{(v)}$ is central in
$\Z \langle Z_{v},x \rangle $ and is independent
of the order of the conjugates in the product
expression. In particular, $b_{(n)}\in Z(\U(\Z
G))$. Let $$B_{(n)}(G)= \left\langle b_{(n)} \mid
b \mbox{ a Bass cyclic of } \Z G\right\rangle.$$

\begin{proposition} \cite[Proposition 2]{JPS} \label{Bn}
If $G$ is a finite nilpotent group of class $n$, then $B_{(n)}(G)$ has finite
index in $Z(\U(\Z G))$.
\end{proposition}

The proof of the previous proposition relies on results of Bass \cite[Lemma~2.2, Lemma~3.6, Theorem~2, Theorem~4]{Bas},
which states that the natural images of the Bass cyclic units in the Whitehead group $K_1(\Z G)$ of $\Z G$ generate a
subgroup of finite index, and on the fact that the torsion-free rank of the abelian groups $K_1(\Z G)$ and $Z(\U(\Z
G))$ are the same. Clearly, because $K_1(\Z G)$ is commutative, the natural image of a $b_{(n)}$ in $K_1(\Z G)$ is
equal with the natural image of some power of $b$ in $K_1(\Z G)$. Hence in $K_1(\Z G)$ the group generated by the Bass
cyclic units contains the group generated by the $b_{(n)}$'s as a subgroup of finite index. So, indeed, $B_{(n)}(G)$ is
of finite index in $Z(\U(\Z G))$.

\vspace{2mm} If $G$ is a finite group and $e$ is a primitive central idempotent
of the rational group algebra $\Q G$, then the simple algebra $\Q Ge$ is
identified with $M_{n}(D)$, a matrix algebra over a division algebra $D$. As in
\cite{JL1}, an exceptional component of $\Q G$ is a non-commutative division
algebra other than a totally definite quaternion algebra, or a two-by-two
matrix algebra over either the rationals, a quadratic imaginary extension of
the rationals or a non-commutative division algebra.

Let $\Or$ be an order in $D$ and denote by $\GL_{n}(\Or)$ the group of
invertible matrices in $M_{n}(\Or)$ and by $\SL_{n}(\Or)$ its subgroup
consisting of matrices of reduced norm $1$. For an ideal $Q$ of $\Or$ we denote
by $E(Q)$ the subgroup of $\SL_{n}(\Or)$ generated by all $Q$-elementary
matrices, that is $E(Q)=\langle I+ qE_{ij} \mid  q\in Q, 1\leq i,j\leq n, i\neq
j, E_{ij} \,\text{\rm a matrix unit} \rangle$. We recall the following
celebrated theorem (see for instance \cite{RS1} or \cite[Theorem 2.2]{JL1}).

\begin{theorem}[Bass-Milnor-Serre-Vaserstein] \label{congruence} If $n\geq 3$ then
$[\SL_{n}(\Or):E(Q)]<\infty$. If $n=2$ and $D$ is an algebraic number field
which is not rational or imaginary quadratic, then
$[\SL_{2}(\Or):E(Q)]<\infty$.
\end{theorem}

Because of the description of the matrix units from
Corollary~\ref{matrixunits}, we can now show that, in case $G$ is nilpotent and
does not have exceptional simple components, $\U(\Z G)$ has a subgroup of
finite index that is generated by three nilpotent groups, one of which is a
central subgroup contained in the group generated by the Bass cyclic units and
the others are generated by the units of the form $1+E_{tt}gE_{t't'}$, with
$g\in G$ and $t,t'\in T_e$, with $T_e$ as in Notation~\ref{Tx}.

\begin{theorem} \label{finiteindex}
Let $G$ be a finite nilpotent group of class $n$ such that $\Q G$ has no
exceptional components. For every primitive central idempotent $e=e(G,H,K)$,
with $(H,K)$ a strong Shoda pair of $G$, let $N=N_G(K)$,
$\varepsilon=\varepsilon(H,K)$ and let $T_e$ and $\beta_e$ be as in
Notation~\ref{Tx}. Fix an order $<$ in $T_e$. Let $H/K=\langle
\overline{a}\rangle$ and let $l$ be the least integer such that
$a^l\varepsilon$ is central in $\Q N\varepsilon$. Then the following two groups
are nilpotent subgroups of $\U(\Z G)$:
    $$V^+_e=\left\langle 1 + |G| t^{-1} \beta_e a^j t' \mid
    a^j\in \langle a^l\rangle, t,t'\in T_e, t'>t\right\rangle,$$
    $$V^-_e=\left\langle 1 + |G| t^{-1} \beta_e a^j t' \mid  a^j\in \langle
    a^l\rangle, t,t'\in T_e, t'<t\right\rangle.$$
Hence
\begin{eqnarray*}V^+=\left\langle V^+_e \mid  e=e(G,H,K) \text{ \rm a primitive central
idempotent of }\Q G \right\rangle =\prod_e V^+_e \text{ \rm and } \\
V^-=\left\langle V^-_e \mid  e=e(G,H,K) \text{ \rm a primitive central
idempotent of }\Q G \right\rangle =\prod_e V^-_e \quad\quad\,
\end{eqnarray*} are nilpotent subgroups of $\U (\Z G)$. Furthermore, if
$B_{(n)}(G)=\langle b_{(n)} \mid  b \text{ a Bass cyclic unit of } G \rangle$, where $n$ is the nilpotency class of
$G$, then the group
 $$\left\langle B_{(n)}(G), V^+, V^-\right\rangle$$
is of finite index in $\U(\Z G)$.
\end{theorem}

\begin{proof}
Recall that the intersection of the unit groups of two orders in a finite dimensional rational algebra are
commensurable and henceforth it is enough to show that $\GEN{B_{(n)},V^+,V^-}$ contains a subgroup of finite index in
the group of units of an order of $(1-e)\Q+\Q Ge$ for a primitive central idempotent $e$ of $\Q G$. So fix such a
primitive central idempotent $e=e(G,H,K)$ of $\Q G$.

The elements of the form $1+ |G| t^{-1}\beta_e a^j t'$,  with $a^j\in \langle a^l \rangle$ and $t,t'\in T_e$, project
trivially to $\Q G(1-e)$ and by Corollary~\ref{matrixunits} they project to an elementary matrix of $M_n(\Or)$, for
some order $\Or$ in the division ring $D$, where $\Q Ge\simeq M_n(D)$. Since also $|G| \beta_e\in \Z G$,  it follows
that $1 +|G| t^{-1} \beta_e a^j t'\in \U(\Z G)$ and by Theorem~\ref{congruence}, for $t\neq t'$, these units generate a
subgroup of finite index in $(1-e)+\SL_n(\Or)$. By Theorem~\ref{Bn}, $B_{(n)}$ has finite index in $Z(\U(\Z G))$ and
therefore it contains a subgroup of finite index in the center of $(1-e)+\GL_n(\Or)$. As the center of $\GL_n(\Or)$
together with $\SL_n(\Or)$ generate a subgroup of finite index of $\GL_n(\Or)$, we conclude that the group
$\GEN{B_{(n)},V^+,V^-}$ contains a subgroup of finite index in  the group of units of an order of $(1-e)\Q+\Q Ge$. Note
that $V_e^+$ and $V_e^-$ correspond to upper and lower triangular matrices respectively and hence they are nilpotent
groups.
\end{proof}

Another application of the construction of the matrix units is that one can
easily obtain free subgroups of $\U(\Z G)$ for $G$ a finite nilpotent group.

\begin{corollary}
Let $G$ be a  finite nilpotent group, $(H,K)$ a
strong Shoda pair of $G$,
$\varepsilon=\varepsilon(H,K)$, $e=e(G,H,K)$ and
let $T_e$ and $\beta_e$ be as in
Notation~\ref{Tx}. If $\Q G e$ is not a division
algebra (see Remark~\ref{simple}), then for every
$t,t'\in T_e$ with $t\neq t'$, $$\left\langle
1+|G|  t^{-1} \beta_e t',  1+|G| t'^{-1} \beta_e
t\right\rangle $$ is a free group of rank~$2$.
\end{corollary}
\begin{proof}
By Corollary~\ref{matrixunits}, we may write
$1+|G| t^{-1} \beta_e t'=1+ |G| E_{tt'}$ and
$1+|G|  t'^{-1} \beta_e t=1+ |G| E_{t't}$. Hence
$\left\langle 1+|G| t^{-1} \beta_e t', 1+|G|
t'^{-1} \beta_e t \right\rangle$ is isomorphic with $\left\langle \left(%
\begin{array}{cc}
  1 & |G| \\
  0 & 1 \\
\end{array}%
\right), \left(%
\begin{array}{cc}
  1 & 0 \\
  |G| & 1 \\
\end{array}%
\right) \right\rangle$. Since $|G|\geq 2$, a well known result of Sanov yields that this group is a free of rank~$2$.
\end{proof}

\begin{remark} \rm
A well known result of Hartley and Pickel \cite{HP} (or see for example
\cite{Seh1}) says that $\U(\Z G)$ contains a free non-abelian subgroup for any
finite non-abelian group $G$ that is not a Hamiltonian $2$-group. Only in 1997,
Sehgal and Marciniak in \cite{MS} gave a concrete construction of such a group.
They showed that if $u_{g,h}=1+(1-g)h\widehat{g}|\langle g\rangle| $ is a non
trivial bicyclic unit then $\left\langle u_{g,h},
u'_{g^{-1}h^{-1}}\right\rangle$ is a free group of rank $2$. More generally, in
\cite{MS2} it is shown that if $c\in\Z G$ satisfies $c^2=0$ and $c\not=0$ then
$\left\langle 1+c^*,1+c\right\rangle$ is free group of rank $2$. For
$c=\sum_{g\in G} z_{g}g$ we denote by $c^{*}=\sum_{g\in G}z_{g} g^{-1}$. So,
$*$ denotes the classical involution on $\Z G$.  The above corollary yields
many more concrete elements that can be substituted for $c$. Since then, as
mentioned in the introduction, there have been several papers on constructing
free subgroups in $\U(\Z G)$ generated by Bass and/or bicyclic units.
\end{remark}

\medskip

\noindent {\bf Acknowledgements.} The second author would like to thank for the warm hospitality during the visit to
Vrije Universiteit Brussel with a postdoctoral grant of Fundaci\'on S\'eneca of Murcia. The authors would like to thank
Capi Corrales for the help with the splitting of quaternion algebras (Lemma~\ref{split quaternion}).

\end{document}